\documentclass[11pt]{article}

\usepackage{bm}
\usepackage{epsfig}
\usepackage{caption}
\usepackage{amsmath,amssymb,amsthm,amsfonts} % Typical maths resource packages

\usepackage{graphics}             % Packages to allow inclusion of graphics
\usepackage{graphicx}
\usepackage{array}
\usepackage{epsf}
\usepackage{color,soul}         % For creating coloured text and background
\usepackage{xcolor}
\usepackage{accents}
\usepackage{empheq}

\newtheorem{theorem}{Theorem}[section]
\newtheorem{lemma}[theorem]{Lemma}
\newtheorem{remark}[theorem]{Remark}

\newtheorem{example}{Example}[section]

\newcommand{\nn}{\nonumber}

\def\refe#1{(\ref{#1})}
\def \bb#1{\setbox0=\hbox{$#1$}
\kern-.025em\copy0\kern-\wd0 \kern.05em\copy0\kern-\wd0
\kern-.025em\raise.0433em\box0}

\newcommand{\D}{\mathrm{div}}

\newcommand{\bH}{\mathbf{H}}
\newcommand{\ee}{\epsilon}
\begin{document}

\date{\today}

\title{\bf Error analysis of mixed finite element methods
  for nonlinear parabolic equations}

\author{
\setcounter{footnote}{0}
Huadong~Gao
\footnote{
School of Mathematics and Statistics, 
Huazhong University of Science and Technology, 
Wuhan 430074, People¡¯s Republic of China.
{\tt huadong@hust.edu.cn}.
The work of the author was supported in part by a grant
from the National Natural Science Foundation of China (NSFC) under
grant No. 11501227.
}
\quad and \quad Weifeng~Qiu
\footnote{ 
Department of Mathematics,
City University of Hong Kong, Kowloon, Hong Kong.
{\tt weifeqiu@cityu.edu.hk}.
The work of the author was supported in part by a grant
from the Research Grants Council of the Hong Kong Special
Administrative Region, China. (Project No. CityU 11302014)  
}
}

\maketitle

\begin{abstract}
  In this paper, we prove a discrete embedding inequality
  for the Raviart--Thomas mixed finite element methods for
  second order elliptic equations,
  which is analogous to the Sobolev embedding inequality
  in the continuous setting.
  Then, by using the proved discrete embedding inequality,
  we provide an optimal error estimate for linearized
  mixed finite element methods for nonlinear parabolic equations.
  Several numerical examples are provided to confirm the theoretical analysis. 

\end{abstract}


\section{Introduction}
\setcounter{equation}{0}

Since the pioneering work of Raviart and Thomas \cite{Raviart_Thomas},
mixed finite element methods (FEMs) have proved to be
a fundamental tool to numerically solve various problems
arising in physics and engineering sciences.
Mixed FEMs have many attractive features over the conventional Lagrange FEMs.
For instance, mixed FEMs conserve mass locally,
which is of crucial importance in numerical methods for flow coupled to transport,
see \cite{Dawson-Sun-Wheeler}.
By introducing $\nabla u$ as an extra variable, mixed methods
can produce accurate flux approximations.
Mixed FEMs are also more robust in the case of low regularity of nonsmooth coefficients.
We refer to the monographs
\cite{Boffi-Brezzi-Fortin,Ern-Guermond} for the general theory of Mixed FEMs.
There are numerous applications of Raviart--Thomas mixed FEMs,
etc., see \cite{Chen_Chen,Garcia1,Garcia2,JohnsonThomee,KimPB, Wu_Allen}
for general linear and nonlinear parabolic equations,
\cite{BrezziMMPSW,Gadau-Jungel,HolstJP} for semiconductor modeling,
\cite{Arbogast-EST,Arbogast-Wheeler} for porous media flow problems.

In this paper, we consider the mixed finite element methods(FEMs) for
nonlinear parabolic equations
\begin{empheq}[left=\empheqlbrace]{align} 
  & \frac{\partial u}{\partial t} - \Delta u + f(u,\nabla u) = 0\, ,
  && \textrm{in $ \Omega \times (0,T]$},
    \label{pde1} \\
  & u(\bm{x},t) = u^0(\bm{x})\,,
  && \textrm{in $ \Omega $},
    \label{pde-initial} \\
  & u = 0\,,
  && \textrm{on $ \partial \Omega $},
    \label{pde-bc}
\end{empheq}
where $f$ is a given function.
We assume $\Omega$ is a bounded Lipschitz polygonal/polyhedral domain
in $\mathbb{R}^d (d=2,3)$.
There have been extensive works on mixed methods for the above equations
\refe{pde1}-\refe{pde-bc}.
We give a brief summary along with some representative
(but certainly not exhaustive) requirement
for the function $f(u,\nabla u)$ (or $f(u)$) in the literature.
\vskip0.1in 

$\bullet$ The function $f(u)$ is a smooth function of $u \in \mathbf{R}$,
i.e., there exists a bound $B_1$ such that
\[
|f^{(i)} | \le B_1, \quad i = 0, 1, 2, \ldots
\]
see page 270 of \cite{KimPB}.

$\bullet$  
There exists a bound $k_1$ such that, for every  $u \in \mathbf{R}$
\[
|f| \le k_1, \quad
\left|\frac{\partial f}{\partial u}\right| \le k_1
\]
See page 131 of \cite{Garcia1} and page 150 of \cite{Garcia2}.

$\bullet$  
The function $f(\bm{x,},u): \overline{\Omega} \times \mathbf{R} \rightarrow \mathbf{R}$
is a triple continuously differential function with bounded derivatives
up to the third order.
See page 205 of \cite{Chen_Chen},
page 410 of \cite{Chen_Liu_Liu} and
page 195 of \cite{Chen_Huang_Yu}.

$\bullet$  
The reaction term $f$ is twice continuously differential on $\Omega$
with bounded derivatives up to the second order.
See page 322 of \cite{Wu_Allen}.

Let us examine why all the above works need these strong assumptions on $f$.
It is not difficult to deduce that, under the above requirements on $f$,
the error between $f(u)$ and $f(u_h)$ can be bounded by
\[
\|f(u) - f(u_h) \|_{L^2} \le \|f'(\xi_{u_h,u})(u - u_h) \|_{L^2}
\le \max_{x \in \mathbf{R}}|f'(x)| \|(u - u_h) \|_{L^2}  \le C \|u - u_h\|_{L^2}.
\]
where $C$ is a constant independent of $u$, $u_h$ and the mesh size $h$.
Due to these severe restrictions on $f$ in
\cite{Garcia1,Garcia2,Chen_Chen,Chen_Liu_Liu,Chen_Huang_Yu,KimPB,Wu_Allen},
the ``nonlinear'' problem \refe{pde1}-\refe{pde-bc} can be almost reduced to a linear one.
Clearly, these assumptions on $f$ cannot be satisfied in most applications.
For instance, $f(u) = u^3-u$ is
frequently used in phase field problems and nonlinear
Schrodinger equations;
$f(u) = ({\bf{b}}\cdot \nabla u)u$ where ${\bf{b}} = [1,1,1]^T$
appears in the viscous Burgers' equation.
In these two cases, $f$ does not satisfies the conditions in
\cite{Chen_Chen,Chen_Liu_Liu,Chen_Huang_Yu,Garcia1,Garcia2,KimPB,Wu_Allen}.
Thus, all previous results are not applicable.
To eliminate the strong assumptions of $f$ and also control the nonlinear term
$f(u, \nabla u)$, one must derive
a uniform boundedness of $u_h$ in certain strong norms.
If conventional Lagrange elements are used, one popular linearized FEM for
the equation \refe{pde1} is to seek $u_h^n \in P_h^r \subset H^1(\Omega)$ such that
\begin{align} 
  & \Big(\frac{u_h^n-u_h^{n-1}}{\tau} , v_h \Big)
  + (\nabla u_h^n , \nabla v_h) + (f(u_h^{n-1},\nabla u_h^n),v_h) = 0\, ,
  && \forall v_h \in P_h^r \subset H^1(\Omega)\,.
    \label{general-lagrange-fem}
\end{align}
And it is easy to see that  $u_h \in H^1(\Omega)$ satisfies
\begin{empheq}[left=\empheqlbrace]{align}
  &  \|u_h\|_{L^p}  \le  C\|\nabla u_h\|_{L^2}  \,,
  && \textrm{for $1 \le p < + \infty$},
  \quad \textrm{in two-dimensional space}
  \nn \\
  &  \|u_h\|_{L^p}  \le  C\|\nabla u_h\|_{L^2}  \,,
  && \textrm{for $1 \le p \le 6$},
  \quad \textrm{in three-dimensional space},
  \nn
\end{empheq}
where $\|\nabla u_h\|_{L^2}$ naturally arise in the discretization of the diffusion
terms. By using this technique,
unconditionally optimal error estimates of conventional Lagrange FEMs
were established in \cite{GLS1,Hou-Li-Sun}
for several nonlinear parabolic problems.
On the contrary, if mixed FEMs are used for spatial discretizations,
a linearized mixed FEM is to seek
$(\bm{\sigma}_h^n,u_h^n) \in \bH_h^r \times V_h^r \subset \bH(\D) \times L^2$,
such that
\begin{empheq}[left=\empheqlbrace]{align} 
  & (\bm{\sigma}_h^n,\bm{\chi}_h) + (u_h^n, \D \bm{\chi}_h) =0\, ,
  && \forall \bm{\chi}_h \in \bH_h^r \,,
    \label{general-mixed-fem1} \\
  & \left(\frac{u_h^n-u_h^{n-1}}{\tau} , v_h \right)
    - (\D \bm{\sigma}_h^n,v_h) + (f(u_h^{n-1},\bm{\sigma}_h^{n-1}),v_h) = 0\, ,
  && \forall v_h \in V_h^r\,,
    \label{general-mixed-fem2}
\end{empheq}
where $\bH_h^r \times V_h^r$ are the Raviart--Thomas mixed finite element spaces,
see the definition in section \ref{preparation}.
It is easy to see that
for the mixed scheme \refe{general-mixed-fem1}-\refe{general-mixed-fem2}
one can only derive 
\[
\frac{\|u_h^n\|^2 - \|u_h^n\|^2}{2 \tau}
+ \frac{1}{2\tau} \|u_h^n - u_h^{n-1}\|
+ \|\bm{\sigma}_h\|_{L^2}^2 +
(f(u_h^{n-1},\bm{\sigma}_h^{n-1}), u_h) = 0.
\]
Therefore, we must use $\|\bm{\sigma}_h\|_{L^2}^2$ 
to control the nonlinear term $f(u_h^{n-1},\bm{\sigma}_h^{n-1})$.
More precisely, we shall establish embedding inequalities between
$u_h^n$ and $\bm{\sigma}_h^n$
\begin{empheq}[left=\empheqlbrace]{align}
  &  \|{u}_h^n\|_{L^p}  \le  C \|\bm{\sigma}_h^n\|_{L^2}  \,,
  \quad \textrm{for $1 \le p < \infty$},
  \quad \textrm{in two-dimensional space},
  \nn \\
  & \|{u}_h^n\|_{L^p}  \le  C\|\bm{\sigma}_h^n\|_{L^2}  \,,
  \quad \textrm{for $1 \le p \le 6$},
  \quad \textrm{in three-dimensional space},
  \nn
\end{empheq}
Although the above results seems rather reasonable,
to the best our knowledge, 
such a relationship for $(\bm{\sigma}_h , u_h)$
is unavailable in the literature.
In this paper, we provide a rigorous proof for the discrete Sobolev embedding
inequalities for the Raviart--Thomas mixed FEMs.
A key step in our proof is to introduce
a new norm $\|u_h\|_{DG}$, which can be viewed as the broken
$H^1$ norm of $u_h$. Then, we analyze $\|u_h\|_{DG}$ and $\|\bm{\sigma}_h\|$
carefully to derive the desired results.

The rest of this paper is organized as follows. 
In section \ref{preparation}, we provide some notations and lemmas for later use.
In section \ref{sobolev-inequality},
we prove the discrete Sobolev embedding inequalities
associated to the Raviart--Thomas mixed FEMs.
In section \ref{error-estimate}, we provide an
optimal error estimate for the linearized mixed FEMs
by using the discrete Sobolev embedding inequalities.
Numerical examples in both two- and three-dimensional spaces
are given in section \ref{sec-numer} to
confirm the error estimates and show the efficiency of linearized mixed FEMs.

\section{Preliminaries}
\label{preparation}
\setcounter{equation}{0}

Let $W^{k,p}(\Omega)$ be the Sobolev space defined on $\Omega$, 
and by conventional notations, $H^{k}(\Omega) := W^{k,2}(\Omega)$.
To introduce the mixed formulation,
we denote 
\begin{align}
{\bH}(\D;\Omega) = \left\{
\bm{u} \, \big| \, \bm{u} \in \mathbf{L}^{2}(\Omega),
\D \, \bm{u}  \in L^2(\Omega) \right\} \,
\textrm{with $\|\bm{u}\|_{{\bH}(\D)} = \left(\|\bm{u}\|_{{L}^2}^2 
+ \|\D \, \bm{u}\|_{{L}^2}^2  \right)^{\frac{1}{2}} $}
\nn
\end{align}
and its dual space $\accentset{\circ}{\bH}(\D)'$
with norm 
\begin{align}
\left \| \bm{v} \right \|_{\accentset{\circ}{\bH}(\D)'}
:= \sup_{\bm{w} \in \accentset{\circ}{\bH}(\D)}
\frac{(\bm{v} \, , \, \bm{w})}
{\quad \left \| \bm{w} \right \|_{{\bH}(\D)}}\,.
\nn
\end{align}

Let $\mathcal{T}_h = \{K\}$ be a regular mesh partition of $\Omega$
and denote the mesh size $h = \max_{K} \{\mathrm{diam} K \}$.
By $\mathcal{F}_h$ we denote all the $(d-1)$-dimensional faces of the partition
$\mathcal{T}_h$.
We define the Raviart--Thomas finite element spaces by 
\begin{empheq}[left=\empheqlbrace]{align} 
  & \bH_h^r(\Omega)  := \{\bm{q} \in \bH(\D;\Omega) \quad : \quad
  \bm{q}|_{K} \in [P_r(K)]^d + \bm{x} P_r(K), \quad \forall K \in \mathcal{T}_h\} \,,
  \nn \\
  & V_h^r(\Omega)  :=  \{u \in L^{2}(\Omega) \quad : \quad
  u|_{K} \in P_r(K), \quad \forall K \in \mathcal{T}_h\} \,,
  \nn 
\end{empheq}
where $P_r(K)$ is the space of polynomials of degree $r$ or less defined on K.
It is well-known that $\bH_h^r(\Omega) \times V_h^r(\Omega)$ is a stable 
finite element pair for solving the second order elliptic problems,
see \cite{Boffi-Brezzi-Fortin,JohnsonThomee,Nedelec,Raviart_Thomas}.
Let ${\left \{ t_{n} \right \}}_{n=0}^{N}$ be a uniform partition
in the time direction with the step size $\tau = \frac{T}{N}$. 
For a sequence of functions $\{u^{n}\}$ defined on $\Omega$, we denote
the backward Euler discretization operator
\begin{eqnarray}
{D_{\tau}} u^{n} = \frac{u^{n}-u^{n-1}}{\tau}.
\nn
\end{eqnarray}

In the rest part of this paper, for simplicity of notation
we denote by $C$ a generic positive
constant and $\ee$ a generic small positive constant, 
which are independent of $n$, $h$ and $\tau$. 
We present the Gagliardo--Nirenberg and discrete Gronwall's inequalities
in the following lemmas which will be frequently used in our proofs.
\begin{lemma}
\label{GN} 
{\it ( Gagliardo--Nirenberg inequality \cite{Nirenberg}): 
Let $u$ be a function defined on $\Omega$ in $\mathbb{R}^d$ 
and $\partial ^{s} u$ be any partial derivative of $u$ of order $s$, then
\begin{equation}
\|\partial ^{j} u\|_{L^p} 
\le C \|\partial^{m} u\|_{L^r}^{a} \, \|u\|_{L^q}^{1-a}
+ C \|u\|_{L^q},
\nn
\end{equation}
for $0 \le j < m$ and $\frac{j}{m} \le a \le 1$ with
\[
\frac{1}{p} = 
\frac{j}{d} + a \left( \frac{1}{r} - \frac{m}{d}\right)
+(1-a) \frac{1}{q} \, ,
\]
except $1 < r < \infty$ and $m-j-\frac{n}{r}$ 
is a non-negative integer,
in which case the above estimate 
holds only for $\frac{j}{m} \le a < 1$.
}
\end{lemma} 

%------------------------------
%---------- Discrete Gronwall's inequality
%------------------------------
\begin{lemma}
\label{gronwall}
Discrete Gronwall's inequality {\cite{Heywood_Rannacher}} : 
Let $\tau$, $B$ and $a_{k}$, $b_{k}$, $c_{k}$, $\gamma_{k}$, 
for integers $k \geq 0$, be non-negative numbers such that
\[
a_{n} + \tau \sum_{k=0}^{n} b_{k} 
\leq \tau \sum_{k=0}^{n} \gamma_{k} a_{k} + 
\tau \sum_{k=0}^{n} c_{k} + B \, , \quad \mathrm{for } \quad n \geq 0 \, ,
\]
suppose that $\tau \gamma_{k} < 1$, for all $k$,
and set $\sigma_{k}=(1-\tau \gamma_{k})^{-1}$. Then
\[
a_{n} + \tau \sum_{k=0}^{n} b_{k} 
\leq  \exp(\tau \sum_{k=0}^{n} \gamma_{k} \sigma_{k}) (\tau \sum_{k=0}^{n} c_{k} + B) \, , 
\quad \mathrm{for } \quad n \geq 0 \, .
\]
\end{lemma}

\section{Discrete Sobolev embedding inequalities of mixed FEMs for the Poisson problem}
\label{sobolev-inequality}
\setcounter{equation}{0}

We consider in this section the model problem
\begin{empheq}[left=\empheqlbrace]{align} 
  & -\Delta u = f\, ,
  &&  \textrm{for $x \in \Omega$}\,,
  \nn\\
  & u = 0\, ,
  &&  \textrm{for $x \in \partial\Omega$}\,.
    \nn
\end{empheq}
The standard Raviart--Thomas mixed FEM for the above model problem
is to seek
$(\bm{\sigma}_h, u_h) \in \bH_h^r(\Omega) \times V_h^r(\Omega) $ such that
\begin{empheq}[left=\empheqlbrace]{align} 
  & (\bm{\sigma}_h,\bm{\chi}_h)
  + (u_h, \D \bm{\chi}_h) =0\, ,
  && \forall \bm{\chi}_h \in \bH_h^r(\Omega) \,,
  \label{poisson-fem1} \\
  &  (\D \bm{\sigma}_h,v_h) =
  - (f,v_h)\, ,
  && \forall v_h \in V_h^r(\Omega) \,.
  \label{poisson-fem2}
\end{empheq}
Error analyses for the above mixed methods \refe{poisson-fem1}-\refe{poisson-fem2}
can be found in \cite{Boffi-Brezzi-Fortin,Nedelec,Raviart_Thomas}
and references therein.
The mixed methods computes the original unknown $u_h$ and the flux $\bm{\sigma}_h$ 
simultaneously.
On the contrary, to obtain the flux $\nabla u$, 
conventional Lagrange FEMs need to use certain numerical differentiation,
which may lead to a loss in accuracy.
If we still denote by $u_h$ the conventional Lagrange FEM solution
to the above Poisson problem,
then $u_h$ satisfies the following Sobolev inequalities
\begin{empheq}[left=\empheqlbrace]{align}
  &  \|{u}_h\|_{L^p}  \le  C \|\nabla u_h\|_{L^2}  \,,
  \quad \textrm{for $1 \le p < \infty$},
  \quad \textrm{in two-dimensional space},
  \label{embedding-p1-2d} \\[3pt]
  & \|{u}_h\|_{L^p}  \le  C \|\nabla u_h\|_{L^2}  \,,
  \quad \textrm{for $1 \le p \le 6$},
  \quad \textrm{in three-dimensional space},
  \label{embedding-p1-3d}
\end{empheq}
which inherits the $H^1$ conforming nature.
As $\bm{\sigma}_h$ numerically converges to $\nabla u$,
one may ask whether similar Sobolev embedding inequalities
hold for the mixed FEM solutions $(\bm{\sigma}_h, u_h)$.
In this section, we give an affirmative answer to this question
and provide a proof.

The main idea used in the proof is to investigate the relationship between
mixed FEMs and the discontinuous Galerkin FEMs.
The reasons are twofold: Firstly, there have been powerful tools developed for the
discontinuous Garkerkin methods, see \cite{BuffaOrtner,DiPietro-Ern};
Secondly, the numerical solution $u_h$ is in the discontinuous finite element space
but not in $H^1(\Omega)$.
Following \cite{DiPietro-Ern}, we define the $\|\cdot\|_{DG}$ norm of $u_h$ by
\begin{equation}
  \|u_h \|_{DG}^2 := \sum_{K \in \mathcal{T}_h} \int_{K} |\nabla u_h|^2 \textrm{d}x
  + \sum_{F \in \mathcal{F}_h} \frac{1}{h_F} \int_{F} |[\![{u}_h]\!]|^2 \textrm{d}x,
\end{equation}
where $h_F$ denotes the size of the face $F$.
For two adjacent elements $K$ and $K'$ sharing the same face $F$,
the jump of a function $u_h \in V_h^r(\Omega)$ across $F$ is defined by
\[
[\![{u}_h]\!] = u_h|_{\partial K \cap F}  -  u_h|_{\partial K' \cap F} .
\]
In case of $F \in \partial K$ lying on $\partial \Omega$, we define
$[\![{u}_h]\!] := u_h|_{\partial K \cap F}$.
The main results obtained in \cite{DiPietro-Ern} are the following
discrete Sobolev embedding inequalities
\begin{empheq}[left=\empheqlbrace]{align}
  &  \|{u}_h\|_{L^p}  \le  C \|{u}_h\|_{DG}  \,,
  \quad \textrm{for $1 \le p < \infty$},
  \quad \textrm{in two-dimensional space},
  \label{embedding-results-2d-DG} \\[3pt]
  & \|{u}_h\|_{L^p}  \le  C\|{u}_h\|_{DG}  \,,
  \quad \textrm{for $1 \le p \le 6$},
  \quad \textrm{in three-dimensional space},
  \label{embedding-results-3d-DG}
\end{empheq}
where $C$ is a constant depending upon the domain $\Omega$,
$r$ and $p$ only.
We will show in Theorem \ref{discrete-embedding-rt} that
$\|{u}_h\|_{DG}$ is bounded by $\|\bm{\sigma}_h\|_{L^2}$.
Before going further,
we consider the following projection problem for the Raviart--Thomas element,
which will be used in the later proof.
This lemma first appears in \cite{EggerS}.
Here we provide a complete proof with details.

\begin{lemma}
  \label{projection-stability}
For each element $K \in \mathcal{T}_h$,
given $\bm{p} \in [{L}^2(K)]^d$, $q_i \in L^2(F_i)$ where $\{F_i\} \in \partial K$,
there exists a unique ${\bm{\zeta}}_h \in RT_r(K)$ such that
\begin{empheq}[left=\empheqlbrace]{align} 
  & \int_{K}({\bm{\zeta}}_h - \bm{p}) \cdot \bm{\omega}_h \textrm{d}x= 0\, ,
  && \forall \bm{\omega}_h \in [P_{r-1}(K)]^d \,,
  \label{interpolate-fem1} \\
  &  \int_{F_i}  ({\bm{\zeta}}_h \cdot \mathbf{n}_{F_i} - q_i) \, \mu_h
  \textrm{d}x = 0\, ,
  && \forall \mu_h \in P_r(F_i)  \,,
  \label{interpolate-fem2}
\end{empheq}
where $RT_r(K)$ is the restriction of the
Raviart--Thomas element space $\bH_h^r(\Omega)$ on $K$.
More importantly, the following stability holds
\begin{equation}
  \| {\bm{\zeta}}_h \|_{L^2(K)}^2 \le
  C\left(\| \bm{p} \|_{L^2(K)}^2
  + \sum_{F_i \in \partial K} h \| q_i \|_{L^2(F_i)}^2 \right) \,,
  \label{interpolate-stability}
\end{equation}
where $C$ is independent of $h$ and $K$.
\end{lemma}
\noindent {\bf{ \em Proof:}}
The equation numbers in \refe{interpolate-fem1}-\refe{interpolate-fem2} 
satisfy 
\[
\underbrace{d \binom{d+r-1}{r-1}}_{\textrm{eqn. \refe{interpolate-fem1}}}
\quad + \quad 
\underbrace{(d+1) \binom{d+r-1}{r}}_{\textrm{eqn. \refe{interpolate-fem2}}}
\quad = \quad
\underbrace{\frac{(d+r+1)\, (d+r-1)!}{(d-1)! r!}}_{\mathrm{dim} \{RT_r(K)\}} \,\,,
\]
which immediately yields the the existence and uniqueness of the projection.  
We prove \refe{interpolate-stability} by a scaling argument.
To do so, let $\widehat{K}$ be the reference element,
which can be a simplex in $\mathbb{R}^d$.
Given $\widehat{\bm{p}} \in [{L}^2(\widehat{K})]^d$,
$\widehat{q}_i \in L^2(\widehat{F}_i)$ where $\{\widehat{F}_i\} \in \partial \widehat{K}$,
the corresponding projection on  $\widehat{K}$ is to
seek ${\widehat{\bm{\zeta}}}_h \in RT_r(\widehat{K})$ such that
\begin{empheq}[left=\empheqlbrace]{align} 
  & \int_{\widehat{K}}(\widehat{{\bm{\zeta}}}_h - \widehat{\bm{p}}) \cdot
  \widehat{\bm{\omega}}_h \textrm{d}\hat{x}= 0\, ,
  && \forall \widehat{\bm{\omega}}_h \in [P_{r-1}(\widehat{K})]^d \,,
  \label{interpolate-fem1-refer} \\[5pt]
  &  \int_{\widehat{F}_i}  ({\widehat{\bm{\zeta}}}_h \cdot \mathbf{n}_{\widehat{F}_i}
  - \widehat{q}_i) \, \widehat{\mu}_h
  \textrm{d}\hat{x} = 0\, ,
  && \forall \widehat{\mu}_h \in P_r(\widehat{F}_i)\,.
  \label{interpolate-fem2-refer}
\end{empheq}
The existence and uniqueness of $\widehat{{\bm{\zeta}}}_h$ are obvious.
Furthermore, it is easy to derive that
\begin{equation}
  \left\| \widehat{{\bm{\zeta}}}_h \right\|_{L^2(\widehat{K})}^2 \le
  C\Big( \left\| \widehat{\bm{p}} \right\|_{L^2(\widehat{K})}^2
  + \sum_{\widehat{F}_i \in \partial \widehat{K}}
  \left\| \widehat{q}_i \right\|_{L^2(\widehat{F}_i)}^2 \Big)
  \label{interpolate-stability-refer}
\end{equation}
where $C$ depends upon $\widehat{K}$ and $r$ only.
To build a connection between $\widehat{K}$ and ${K}$,
we define the affine mapping $T_K : \mathbb{R}^d \rightarrow \mathbb{R}^d$ by
\begin{align}
  T_K := B_K \widehat{x} + b_k, \quad \forall \widehat{x} \in \mathbb{R}^d.
  \label{geometry-transform}
\end{align}
We also note that for $\{\widehat{F}_i\} \in \partial \widehat{K}$,
$T_K(\widehat{F}_i) = {F}_i$.
Then, for a scalar function $\phi \in H^1(K)$, we define
\begin{align}
  \widehat{\phi} = \phi \circ T_K
  \label{standard-transform}
\end{align}
For $\bm{\zeta} \in [H^1(K)]^d$, we shall introduce the Piola transformation
\begin{equation}
  \widehat{\bm{\zeta}} = |\mathrm{det} B_K| \, B_K^{-1}   \, \bm{\zeta} \circ T_K .
  \label{piola-transform}
\end{equation}
With the above transformations \refe{geometry-transform}-\refe{piola-transform},
the projection \refe{interpolate-fem1}-\refe{interpolate-fem2}
can be rewritten by
\begin{empheq}[left=\empheqlbrace]{align} 
  & \int_{\widehat{K}}\left(\widehat{{\bm{\zeta}}}_h - \widehat{\bm{p}}\right) \cdot
  \widehat{\bm{\omega}}_h \textrm{d}\hat{x}= 0\, ,
  && \forall \widehat{\bm{\omega}}_h \in [P_{r-1}(\widehat{K})]^d \,,
  \label{interpolate-fem1-trans} \\[3pt]
  &  \int_{\widehat{F}_i}  \left({\widehat{\bm{\zeta}}}_h \cdot \mathbf{n}_{\widehat{F}_i}
  - \frac{|F_i|}{|\widehat{F}_i|} \, \widehat{q}_i \right) \, \widehat{\mu}_h
  \textrm{d}\hat{x} = 0\, ,
  && \forall \widehat{\mu}_h \in P_r(\widehat{F}_i)\,.
  \label{interpolate-fem2-trans}
\end{empheq}
where the Piola transformation \refe{piola-transform} is used for the
vector functions $\bm{\zeta}_h$, $\bm{p}$ and $\bm{\omega}_h$,
and affine transformation \refe{standard-transform} is used
for $\{{q}_i\}$ and $\mu_h$, respectively.
Finally, we use the scaling argument to derive that
\begin{eqnarray}
  \| {\bm{\zeta}}_h \|_{L^2(K)}
  &  \le &
  \|B_K\| \, {|\mathrm{det} B_K|}^{-\frac{1}{2}}
  \| {\widehat{\bm{\zeta}}}_h \|_{L^2(\widehat{K})}
  \nn \\
  &  \le &
  C\|B_K\| \, {|\mathrm{det} B_K|}^{-\frac{1}{2}}
  \left( \left\|\widehat{\bm{p}} \right\|_{L^2(\widehat{K})}
  + \sum_{\widehat{F}_i \in \partial \widehat{K}}\frac{|F_i|}{|\widehat{F}_i|}
  \left \|\widehat{q}_i\right\|_{L^2(\widehat{F}_i)}
  \right)
  \nn \\
  &  \le &
  C\|B_K\| \, {|\mathrm{det} B_K|}^{-\frac{1}{2}}
  \left( |\mathrm{det} B_K|^{\frac{1}{2}} \, \|B_K^{-1}\|
  \left\|{\bm{p}} \right\|_{L^2({K})}
  + \sum_{{F}_i \in \partial {K}}\left(\frac{|F_i|}{|\widehat{F}_i|}\right)^{1/2}
  \left \|{q}_i\right\|_{L^2({F}_i)}
  \right)  
  \nn \\
  &  \le &
  C\|B_K\| \, \|B_K^{-1}\|  \left\|{\bm{p}} \right\|_{L^2({K})}
  + C \sum_{{F}_i \in \partial {K}}
  \|B_K\| \, {|\mathrm{det} B_K|}^{-\frac{1}{2}}
  \left(\frac{|F_i|}{|\widehat{F}_i|}\right)^{1/2}
  \left \|{q}_i\right\|_{L^2({F}_i)}
  \nn \\
  &  \le & C \left(
  \left\|{\bm{p}} \right\|_{L^2({K})}
  + h^{\frac{1}{2}} \sum_{{F}_i \in \partial {K}}\left \|{q}_i\right\|_{L^2({F}_i)}
  \right),
  \nn
\end{eqnarray}
where we have used estimate \refe{interpolate-stability-refer}
and noted the fact that $\mathcal{T}_h$ is shape regular( $h$ and $h_F$ are equivalent).
By squaring the above inequality,
we get the desired estimate \refe{interpolate-stability} directly. \qed

Now, we are ready to prove the discrete Sobolev embedding inequalities
for the mixed FEM solutions $(\bm{\sigma}_h,u_h)$,
which play a key role in the unconditionally optimal error analysis in
the next section \ref{error-estimate}.

\begin{theorem}
  \label{discrete-embedding-rt}
For any given $u_h \in V_h^r(\Omega)$,
if there exists a ${\bm{\sigma}}_h \in \bH_h^r(\Omega)$ such that
\begin{eqnarray}
 ({\bm{\sigma}}_h,\bm{\chi}_h) + ({u}_h, \D \bm{\chi}_h) =0\, ,
  \qquad \forall \bm{\chi}_h \in \bH_h^r(\Omega) \,,
  \label{property-dis-gradient}
\end{eqnarray}
then the following discrete Sobolev embedding inequality holds
\begin{empheq}[left=\empheqlbrace]{align}
  &   \|{u}_h\|_{L^p}  \le  C\|{\bm{\sigma}}_h\|_{L^2}  \,,
  &&  \quad \textrm{for $1 \le p < \infty$},
  && \textrm{in two dimensional space},
  \label{embedding-results-2d} \\[3pt]
  &   \|{u}_h\|_{L^p}  \le  C\|{\bm{\sigma}}_h\|_{L^2}  \,,
  &&  \quad \textrm{for $1 \le p \le 6$},
  && \textrm{in three dimensional space},
  \label{embedding-results-3d}
\end{empheq}
where $C$ is a constant only depending upon the domain, $r$ and $p$.
\end{theorem}

\noindent {\bf{ \em Proof:}}
By using integration by parts for the equation \refe{property-dis-gradient},
we have
\begin{eqnarray}
  0 & = & ({\bm{\sigma}}^n,\bm{\chi}_h) + ({u}_h, \D \bm{\chi}_h)
  \nn \\
  & = &
  ({\bm{\sigma}}^n,\bm{\chi}_h)
  + \sum_{K \in \mathcal{T}_h} \int_{K} {u}_h  \D \bm{\chi}_h \mathrm{d}x
  \nn \\
  & = &
  ({\bm{\sigma}}^n,\bm{\chi}_h)
  + \sum_{K \in \mathcal{T}_h} \left(
  -\int_{K} \nabla {u}_h \cdot \bm{\chi}_h  \mathrm{d}x
  +\int_{\partial K} {u}_h \, \bm{\chi}_h \cdot \mathbf{n}_K \, \mathrm{d}s
  \right)
  \nn \\
  & = &
  ({\bm{\sigma}}^n,\bm{\chi}_h)
  -\sum_{K \in \mathcal{T}_h} 
  \int_{K} \nabla {u}_h \cdot \bm{\chi}_h  \mathrm{d}x
  +\sum_{F \in \mathcal{F}_h} 
  \int_{F} [\![{u}_h]\!] \, \bm{\chi}_h \cdot \mathbf{n}_F \mathrm{d}x \,.
  \label{rt-dg-integration-byparts}
\end{eqnarray}
Then, on each element $K$, we require $\bm{\chi}_h$ to be the projection, such that
\begin{empheq}[left=\empheqlbrace]{align} 
  & \int_{K}({\bm{\chi}}_h - \nabla u_h) \cdot \bm{\omega}_h \textrm{d}x= 0\, ,
  && \forall \bm{\omega}_h \in [P_{r-1}(K)]^d \,,
  \nn \\
  &  \int_{F_i}  ({\bm{\chi}}_h \cdot \mathbf{n}_{F_i} + \frac{1}{h_{F_i}}[\![{u}_h]\!]) \, \mu_h
  \textrm{d}x = 0\, ,
  && \forall \mu_h \in P_r(F_i) \,.
  \nn
\end{empheq}
Lemma \ref{projection-stability} tells that such
a $\bm{\chi}_h \in \bH_h^r(\Omega)$ exists and is unique, which also satisfies
\begin{equation}
  \| {\bm{\chi}}_h \|_{L^2(K)}^2 \le
  C\left(\| \nabla {u}_h \|_{L^2(K)}^2
  + \sum_{F_i \in \partial K} \frac{1}{h} \|[\![{u}_h]\!] \|_{L^2(F_i)}^2 \right)
  \nn
\end{equation}
Substituting this ${\bm{\chi}}_h$ into \refe{rt-dg-integration-byparts} yields
\begin{equation}
  \| u_h \|_{DG}^2 \le (\bm{\sigma}_h,{\bm{\chi}}_h)
  \le \|\bm{\sigma}_h \|_{L^2(\Omega)} \, \|\bm{\chi}_h \|_{L^2(\Omega)}
  \le C \|\bm{\sigma}_h \|_{L^2(\Omega)} \, \| u_h \|_{DG} \,.
  \nn
\end{equation}
The discrete Sobolev embedding inequalities
\refe{embedding-results-2d}-\refe{embedding-results-3d}
follows directly from
\refe{embedding-results-2d-DG}-\refe{embedding-results-3d-DG}
and the last inequality.
The theorem is proved. \qed
\begin{remark}
  In the above proof, we can see that ${\bm{\sigma}}_h$ can be replaced by
  any $\bm{f} \in {\bf{L}}^2(\Omega)$. 
\end{remark}  

\section{Applications in unconditionally optimal error estimates of
nonlinear parabolic equations}
\label{error-estimate}
\setcounter{equation}{0}

\subsection{Linearized mixed FEMs for nonlinear parabolic equations}

In this section, we employ the discrete Sobolev
embedding inequalities to establish an
unconditionally optimal error estimates
of linearized mixed FEMs for nonlinear parabolic equations.
We point out that examples shown in this section do not satisfy the
assumptions assumed in previous works
\cite{Chen_Chen,Chen_Liu_Liu,Chen_Huang_Yu,Garcia1,Garcia2,KimPB,Wu_Allen}.
We present two typical nonlinear equations with different $f(u,\nabla u)$.

\begin{example}
  \label{example-allen-cahn}
The first one is the Allen--Cahn type equation
\begin{empheq}[left=\empheqlbrace]{align} 
  & \frac{\partial u}{\partial t} - \Delta u + u^3-u = 0\, ,
  && \textrm{in $ \Omega \times (0,T]$},
    \label{allen-pde} \\
  & u(\bm{x},t) = u^0(\bm{x})\,,
  && \textrm{in $ \Omega $},
    \label{allen-pde-initial} \\
  & u = 0\,,
  && \textrm{on $ \partial \Omega $},
    \label{allen-pde-bc}
\end{empheq}
\end{example}
Conventional Lagrange FEMs have been widely used to solve the above equation.
However, one might consider using a mixed method for the above Allen--Cahn
equations in the hope of getting a better approximation of the flux
$\nabla u$, which is needed in the case of coupling $\nabla u$ with
Navier--Stokes equations, see \cite{FengHeLiu,Hua_Lin_Liu_Wang}.

\begin{example}
  \label{example-burger}
The second one is the viscous Burgers' equation
\begin{empheq}[left=\empheqlbrace]{align} 
  & \frac{\partial u}{\partial t} - \Delta u + ({\bf{b}} \cdot \nabla u)u = 0\, ,
  && \textrm{in $ \Omega \times (0,T]$},
    \label{burgers-pde} \\
  & u(\bm{x},t) = u^0(\bm{x})\,,
  && \textrm{in $ \Omega $},
    \label{burgers-pde-initial} \\
  & u = 0\,,
  && \textrm{on $ \partial \Omega $}.
    \label{burgers-pde-bc}
\end{empheq}
where ${\bf{b}} = [1,1,1]^T$.
\end{example}

Here we combines the nonlinear terms in the last two examples
and study the following artificial problem
\begin{empheq}[left=\empheqlbrace]{align} 
  & \frac{\partial u}{\partial t} - \Delta u + ({\bf{b}} \cdot \nabla u) \, u
  + u^3 - u = 0\, ,
  && \textrm{in $ \Omega \times (0,T]$},
    \label{artificial-pde} \\
  & u(\bm{x},t) = u^0(\bm{x})\,,
  && \textrm{in $ \Omega $},
    \label{artificial-pde-initial} \\
  & u = 0\,,
  && \textrm{on $ \partial \Omega $}.
    \label{artificial-pde-bc}
\end{empheq}
A linearized mixed FEMs is to look for 
$(\bm{\sigma}_{h}^{n}$, $u_{h}^{n}) \in \bH_h^r(\Omega) \times V_h^r(\Omega)$
such that for $n = 1,2, \ldots$,
\begin{empheq}[left=\empheqlbrace]{align} 
  & (\bm{\sigma}_h^n,\bm{\chi}_h) + (u_h^n, \D \bm{\chi}_h)= 0\, ,
  && \forall \bm{\chi}_h \in \bH_h^r(\Omega) \,,
    \label{artificial-mixed-fem1} \\
  & \left(D_\tau u_h^n, v_h \right)
    - (\D \bm{\sigma}_h^n,v_h) + ({\bf{b}}\cdot \bm{\sigma}_h^{n-1} u_h^{n-1},v_h) 
  &&
  \nn \\
  & \qquad \qquad  \qquad  \qquad  \qquad  \quad   + ((u_h^{n-1})^3-u_h^{n-1}, v_h)= 0\, ,
  && \forall v_h \in V_h^r(\Omega).
  \label{artificial-mixed-fem2}
\end{empheq}
At the initial step, we take $\bm{\sigma}_h^0 = \Pi_h \nabla u^0(\bm{x})$ and
$u_h^0 = \Pi_h u^0(\bm{x})$, where $\Pi_h$  can be the projection defined in
\refe{projection-fem1}-\refe{projection-fem2}.

For error analysis,
we assume that the initial-boundary
value problem (\ref{artificial-pde})-(\ref{artificial-pde-bc}) 
has a unique solution satisfying the regularity condition
\begin{eqnarray}
&& u\in {L^{\infty}(0,T;{H}^{{r+2}})} \, ,
u_{t} \in {L^{\infty}(0,T;{H}^{{r+2}})} \, , 
u_{tt} \in {L^{\infty}(0,T;{L}^2)}.
\label{regularity-artificial}
\end{eqnarray}
We shall remark that the above regularity assumption might be weakened.
In this paper, we emphasize on the unconditional optimal error estimates of
the linearized mixed FEMs \refe{artificial-mixed-fem1}-\refe{artificial-mixed-fem2}.
We state our main results on error analysis in the following theorem.
The proof will be given in the next subsection \ref{proof-parabolic}.

%================
% The main THM for artificial equation
%================
\begin{theorem}
  \label{maintheorem-artificial}
Under the regularity assumption ,
there exist two positive constants
$h_0$ and $\tau_0$ such that when $h<h_0$ and $\tau<\tau_0$, 
the mixed FEM systems \refe{artificial-mixed-fem1}-\refe{artificial-mixed-fem2}
are uniquely solvable and the following error estimates hold
\begin{align}
& \max_{0 \leq n \leq N} \Big(
{\| u_{h}^{n} - u^{n}\|_{L^2}^2}
+ \tau \sum_{m=1}^n{\|\bm{\sigma}_h^m - \bm{\sigma}^m\|_{L^2}^2} \Big)
\leq C_* (\tau^2 + h^{2r+2}) \, ,
\end{align}
where $C_*$ is a positive constant independent of $n$, $h$ and $\tau$. 
\end{theorem}

To prove the above theorem, we need to define a projector
$\Pi_h : (\bH(\D;\Omega),L^2(\Omega)) \rightarrow (\bH_h^r(\Omega), V_h^r(\Omega))$.
Given the exact solution $(\bm{\sigma},u)$ to \refe{artificial-pde}-\refe{artificial-pde-bc}
at any time $t \in (0,T]$,
we seek $(\Pi_h \bm{\sigma},\Pi_h u) \in (\bH_h^r(\Omega),V_h^r(\Omega))$ such that
\begin{empheq}[left=\empheqlbrace]{align} 
  & (\Pi_h\bm{\sigma}_h,\bm{\chi}_h)
  + (\Pi_hu_h, \D \bm{\chi}_h) =0\, ,
  && \forall \bm{\chi}_h \in \bH_h^r(\Omega) \,,
    \label{projection-fem1} \\
 &  (\D  \, \Pi_h\bm{\sigma}_h,v_h) =
    (\D \, \bm{\sigma},v_h)\, ,
  && \forall v_h \in V_h^r(\Omega) \,.
    \label{projection-fem2}
\end{empheq}
We denote the projection error functions by
\begin{equation}
  \theta_{\bm{\sigma}} = \Pi_h \bm{\sigma} - \bm{\sigma},
  \quad
  \theta_{u} = \Pi_h u - u,
\end{equation}
From the classical error estimates for mixed methods
\cite{Boffi-Brezzi-Fortin,Ern-Guermond,Raviart_Thomas},
we have 
\begin{equation}
  \|\theta_{u}\|_{L^2} + \|\theta_{\bm{\sigma}}\|_{L^2} \le Ch^{r+1} \|u\|_{H^{r+1}},
  \quad
  \left\|\frac{\partial \theta_{u}}{\partial t}\right\|_{L^2}
  \le Ch^{r+1} \left\|\frac{\partial u}{\partial t}\right\|_{H^{r+1}}.
  \label{projection-error1}
\end{equation}
Moreover, the following uniform boundedness estimates can be proved by
using an inverse inequality
\begin{equation}
  \left\| \Pi_h \bm{\sigma} \right\|_{L^p}
  + \left\| \Pi_h u \right\|_{L^p} \le C,
  \quad \textrm{for $1 \le p \le 6$}.
  \label{projection-error2}
\end{equation}
With the above projection error estimates, we only need to analyze
the error functions
\begin{equation}
  e_{\bm{\sigma}}^n = \bm{\sigma}_h^n - \Pi_h \bm{\sigma}^n,
  \quad
  e_{u}^n = u^n_h - \Pi_h u^n,
  \quad
  \textrm{for $n=1,2, \ldots, N$}.
\end{equation}
We provide an unconditionally optimal estimates for
$\{(e_{\bm{\sigma}}^n,e_{u}^n)\}_{n=0}^N$ in the next subsection.

\subsection{ Proof of Theorem \ref{maintheorem-artificial}}
\label{proof-parabolic}

\noindent {\bf{ \em Proof:}}
The existence and uniqueness of numerical solutions to the linearized mixed FEMs
\refe{artificial-mixed-fem1}-\refe{artificial-mixed-fem2} follow directly from
that at each time step, the coefficient matrix is invertable.
Here we prove the following inequality for $n = 0$, $\ldots$, $N$
\begin{align}
& {\| e_{u}^{n}\|_{L^2}^2} 
+ \sum_{m=1}^{n} \tau  \| e_{\bm{\sigma}}^m \|_{L^2}^2
\leq \frac{C_*}{2} \left( \tau^{2} + h^{2r+2} \right),
\label{mainresults}
\end{align}
by mathematical induction.
Since 
\begin{align}
& \| e_{u}^{0}\|_{L^2}^2 + \| e_{\bm{\sigma}}^{0}\|_{L^2}^2 = 0,
\nn
\end{align}
(\ref{mainresults}) holds for $n = 0$.
We can assume that (\ref{mainresults}) 
holds for $ n \le k-1$ for some $k\ge 1$.
We shall find a constant $C_*$,
which is independent of $n$, $h$, $\tau$,
such that (\ref{mainresults}) holds for $n \leq k$. 

At time step $t_n$,
by noting the projection \refe{projection-fem1}-\refe{projection-fem2},
the exact solution $(\bm{\sigma},u)$ satisfies
\begin{align} 
  & (\Pi_h{\bm{\sigma}}^n,\bm{\chi}_h) + (\Pi_h{u}^n, \D \bm{\chi}_h) =0\, ,
  && \forall \bm{\chi}_h \in \bH_h^r(\Omega),
  \label{exact-discrete-0} \\[3pt]
  & \left(D_\tau {u}^n, v_h \right)
  - (\D \Pi_h{\bm{\sigma}}^n,v_h)
  = -({\bf{b}}\cdot \bm{\sigma}^{n-1} \, u^{n-1}, v_h)
  && \nn \\
  & \quad \quad \quad \qquad \quad \quad \quad \qquad \quad \quad \quad 
  -((u^{n-1})^3 -u^{n-1},v_h)
  - ( R_{u}^n,v_h)\, ,
  && \forall v_h \in V_h^r(\Omega)\,,
  \label{exact-discrete-1}
\end{align}
where
\begin{align}
  R_{u}^n = D_\tau u^n - \frac{\partial u}{\partial t}\Big|_{t_n}
  + ({\bf{b}}\cdot \nabla u^{n} u^{n} - {\bf{b}}\cdot \nabla u^{n-1} u^{n-1})
  + ((u^{n-1})^3-u^{n-1}) - ((u^{n})^3-u^{n})
  \nn
\end{align}
stands for the truncation error.
Subtracting \refe{exact-discrete-0}-\refe{exact-discrete-1}
from \refe{artificial-mixed-fem1}-\refe{artificial-mixed-fem2},
respectively, we obtain the error equations
\begin{align}
  & (e_{\bm{\sigma}}^n,\bm{\chi}_h) + (e_{u}^n, \D \bm{\chi}_h) =0\, ,
  && \forall \bm{\chi}_h \in \bH_h^r(\Omega),
  \label{allen-error-f1} \\[3pt]
  & \left(D_\tau e_{u}^n, v_h \right)
  - (\D e_{\bm{\sigma}}^n,v_h)
  = ({\bf{b}}\cdot \nabla u^{n-1} u^{n-1}-
  {\bf{b}}\cdot \bm{\sigma}_h^{n-1} u_h^{n-1},v_h)
  && \nn \\
  & \qquad   -(((u^{n-1})^3 -u^{n-1})-((u_h^{n-1})^3 -u_h^{n-1}),v_h)
  -  (D_\tau \theta_{u}^n - R_{u}^n,v_h)\, ,
  && \forall v_h \in V_h^r(\Omega)\,.
  \label{allen-error-f2}
\end{align}
Taking $(\bm{\chi}_h, v_h) = (e_{\bm{\sigma}}^n, e_{u}^n)$ into
the above error equations \refe{allen-error-f1}-\refe{allen-error-f2}
and summing up the results lead to
\begin{align}
  &  \left(D_\tau e_{u}^n , e_{u}^n \right) + \|e_{\bm{\sigma}}^n\|_{L^2}^2
  \nn \\
  &  =
  ({\bf{b}}\cdot \nabla u^{n-1} u^{n-1}-
  {\bf{b}}\cdot \bm{\sigma}_h^{n-1} u_h^{n-1}, e_{u}^n)
  \nn \\
  &  \quad
  + (((u^{n-1})^3 -u^{n-1})-((u_h^{n-1})^3 -u_h^{n-1}), e_{u}^n)
  -  (D_\tau \theta_{u}^n - R_{u}^n, e_{u}^n)
  \label{error-equation}
\end{align}

We estimate the right hand side of \refe{error-equation} term by term.
By the regularity assumption on the exact solution $u$
in \refe{regularity-artificial} and the projection error \refe{projection-error1},
the last term in the right hand side of \refe{error-equation} can be bounded by
\begin{align}
  &  (D_\tau \theta_{u}^n - R_{u}^n, e_{u}^n)
  \le C \|e_{u}^n\|_{L^2}^2 + C \tau^2 + C h^{2r+2}
  \label{term-trunc-1}
\end{align}
By noting that the assumption \refe{mainresults} holds for $ n \le k-1$,
we can derive that if $\tau \le h^{r+1}$
\begin{align}
  & \| e_{u}^n\|_{L^3}
  \le h^{-\frac{1}{2}} \| e_{u}^n\|_{L^2}
  \le h^{-\frac{1}{2}} \sqrt{\frac{C_*}{2} (\tau^{2} + h^{2r+2})}
  \le \sqrt{{C_*}} h^{r+\frac{1}{2}} ,
  \label{assumption-results-1}
\end{align}
where an inverse inequality is used.
If $ h^{r+1}\le \tau$, by using the discrete Sobolev embedding inequality in
Theorem \ref{discrete-embedding-rt}, we have
\begin{align}  
  & \| e_{u}^n\|_{L^3}
  \le C \| e_{\bm{\sigma}}^n\|_{L^2}
  \le C \tau^{-\frac{1}{2}}\sqrt{\tau \| e_{\bm{\sigma}}^n\|_{L^2}^2}
  \le C \tau^{-\frac{1}{2}}\sqrt{\frac{C_*}{2}(\tau^{2} + h^{2r+2})}
  \le C \sqrt{{C_*}} \tau^{\frac{1}{2}} .
  \label{assumption-results-2}
\end{align}
Therefore, for any given $\epsilon$, we have in both cases
\begin{align}
  \| e_{u}^n\|_{L^3}
  \le \max \left\{\sqrt{{C_*}} h^{r+\frac{1}{2}}, C \sqrt{{C_*}} \tau^{\frac{1}{2}}\right\}
  \le \epsilon, \quad
  \textrm{for $n=1,\ldots, k-1$}
  \label{assumption-results}
\end{align}
if we require that $\tau \le \tau_0$ and $h \le h_0$,
where $\tau_0$ and $h$ are two small constant numbers.
Now we estimate the two nonlinear terms in the right hand side of \refe{error-equation}.
The first nonlinear term can be bounded by
\begin{align}
  & ({\bf{b}}\cdot \nabla u^{n-1} u^{n-1}-
  {\bf{b}}\cdot \bm{\sigma}_h^{n-1} u_h^{n-1}, e_{u}^n)
  \nn  \\
  & = -({\bf{b}}\cdot \nabla u^{n-1} (\theta_u^{n-1}+e_u^{n-1}), e_{u}^n)
  -({\bf{b}}\cdot (\theta_{\bm{\sigma}}^{n-1}+e_{\bm{\sigma}}^{n-1}) \Pi_h u^{n-1}, e_{u}^n)
  \nn  \\
  & \quad
  - ({\bf{b}}\cdot (\theta_{\bm{\sigma}}^{n-1}+e_{\bm{\sigma}}^{n-1}) e_u^{n-1}, e_{u}^n)
  \nn  \\
  & \le \|{\bf{b}}\cdot \nabla u^{n-1}\|_{L^\infty}
  \|\theta_u^{n-1}+e_u^{n-1}\|_{L^2} \| e_{u}^n\|_{L^2}
  + \|{\bf{b}}\cdot(\theta_{\bm{\sigma}}^{n-1}+e_{\bm{\sigma}}^{n-1})\|_{L^2}
  \|\Pi_h u^{n-1} \|_{L^6} \|e_{u}^n\|_{L^3}
  \nn  \\
  & \quad 
  + \|{\bf{b}}\cdot (\theta_{\bm{\sigma}}^{n-1}+e_{\bm{\sigma}}^{n-1})\|_{L^2}
  \|e_u^{n-1}\|_{L^3} \| e_{u}^n\|_{L^6}
  \nn  \\
  &
  \le C \| e_{u}^{n-1}\|_{L^2}^2 + C \| e_{u}^n\|_{L^2}^2 + C h^{2r+2}
  + C( \| e_{\bm{\sigma}}^{n-1} \|_{L^2} + h^{r+1})
  \|e_{u}^n\|_{L^2}^{\frac{1}{2}}\|e_{u}^n\|_{L^6}^{\frac{1}{2}}
  \nn  \\
  & \quad  + C( \| e_{\bm{\sigma}}^{n-1} \|_{L^2} + h^{r+1})
  \|e_{u}^{n-1}\|_{L^3}\|e_{u}^n\|_{L^6}
  \nn \\
  & \le C( \| e_{\bm{\sigma}}^{n-1} \|_{L^2} + h^{r+1})
  \|e_{u}^n\|_{L^2}^{\frac{1}{2}}\|e_{\bm{\sigma}}^n\|_{L^2}^{\frac{1}{2}}
  +  C( \| e_{\bm{\sigma}}^{n-1} \|_{L^2} + h^{r+1})
  \|e_{u}^{n-1}\|_{L^3}\|e_{\bm{\sigma}}^n\|_{L^2}
  \nn \\
  & \quad
  + C \| e_{u}^{n-1}\|_{L^2}^2 + C \| e_{u}^n\|_{L^2}^2 + C h^{2r+2}
  \label{nonlinear-term1}
\end{align}
By using Young's inequality, we have
\begin{align}
  & C( \| e_{\bm{\sigma}}^{n-1} \|_{L^2} + h^{r+1})
  \|e_{u}^n\|_{L^2}^{\frac{1}{2}}\|e_{\bm{\sigma}}^n\|_{L^2}^{\frac{1}{2}}
  \le \epsilon (\|e_{\bm{\sigma}}^{n-1}\|_{L^2}^2 + \|e_{\bm{\sigma}}^n\|_{L^2}^2)
  +  \epsilon^{-1} C \|e_{u}^n\|_{L^2}^2 + \epsilon^{-1} C h^{2r+2} \,.
\end{align}
And by using the \refe{assumption-results}, we can derive
\begin{align}
  & C( \| e_{\bm{\sigma}}^{n-1} \|_{L^2} + h^{r+1})
  \|e_{u}^{n-1}\|_{L^3}\|e_{\bm{\sigma}}^n\|_{L^2}
  \le \epsilon ( \| e_{\bm{\sigma}}^{n-1} \|_{L^2}^2
  + \|e_{\bm{\sigma}}^n\|_{L^2}^2 )
  + \epsilon^{-1} C h^{2r+2}
\end{align}
where we shall require $\tau$ and $h$ are smaller than certain constants.
Substituting the last two inequality into \refe{nonlinear-term1} gives
\begin{align}
  & ({\bf{b}}\cdot \nabla u^{n-1} u^{n-1}-
  {\bf{b}}\cdot \bm{\sigma}_h^{n-1} u_h^{n-1}, e_{u}^n)
  \nn \\
  & 
  \le \epsilon ( \| e_{\bm{\sigma}}^{n-1} \|_{L^2}^2
  + \|e_{\bm{\sigma}}^n\|_{L^2}^2 )
  + \epsilon^{-1} C \|e_{u}^n\|_{L^2}^2
  + \epsilon^{-1} C h^{2r+2}
  \label{term-nonlinear-1-final}
\end{align}
Next, the second nonlinear term
in the right hand side of \refe{error-equation} can be bounded by
\begin{align}
  &  (((u^{n-1})^3 -u^{n-1})-((u_h^{n-1})^3 -u_h^{n-1}), e_{u}^n)
  \nn \\
  &  \le ((u^{n-1})^3 - (u_h^{n-1})^3, e_{u}^n)
  + C \|e_{u}^{n-1}\|_{L^2}^2 + C \|e_{u}^n\|_{L^2}^2 + C h^{2r+2}
  \nn \\
  &  \le (-3(u^{n-1})^2(\theta_u^{n-1}+ e_u^{n-1})
  + 3 u^{n-1} (\theta_u^{n-1}+ e_u^{n-1})^2
  -(\theta_u^{n-1}+ e_u^{n-1})^3, e_{u}^n)
  \nn \\
  & \quad  + C \|e_{u}^{n-1}\|_{L^2}^2 + C \|e_{u}^n\|_{L^2}^2 + C h^{2r+2}
  \nn \\
  &  \le  C |((e_u^{n-1})^3, e_{u}^n)| + C|(\theta_u^{n-1}(e_u^{n-1})^2, e_{u}^n)|
  + C \|e_{u}^{n-1}\|_{L^2}^2 + C \|e_{u}^n\|_{L^2}^2 + C h^{2r+2}
  \nn \\
  &  \le  C \|e_u^{n-1}\|_{L^3}\|e_u^{n-1}\|_{L^3}\|e_u^{n-1}\|_{L^6} \|e_{u}^n\|_{L^6}
  + C \|\theta_u^{n-1}\|_{L^3}\|e_u^{n-1}\|_{L^3}\|e_u^{n-1}\|_{L^6} \| e_{u}^n\|_{L^6}
  \nn \\
  & \quad  + C \|e_{u}^{n-1}\|_{L^2}^2 + C \|e_{u}^n\|_{L^2}^2 + C h^{2r+2}
  \nn \\
  &  \le  C \epsilon^2 \|e_{\bm{\sigma}}^{n-1}\|_{L^2} \|e_{\bm{\sigma}}^n\|_{L^2}
  + C \epsilon \|e_{\bm{\sigma}}^{n-1}\|_{L^2} \| e_{\bm{\sigma}}^n\|_{L^2}
  + C \|e_{u}^{n-1}\|_{L^2}^2 + C \|e_{u}^n\|_{L^2}^2 + C h^{2r+2}
  \nn \\
  &  \le \epsilon (\|e_{\bm{\sigma}}^{n-1}\|_{L^2}^2 + \|e_{\bm{\sigma}}^{n}\|_{L^2}^2)
  + C \|e_{u}^{n-1}\|_{L^2}^2 + C \|e_{u}^n\|_{L^2}^2 + C h^{2r+2}
  \label{term-nonlinear-2-final}
\end{align}
where we have used \refe{assumption-results} with requirement $\tau$ and $h$ are smaller
than certain constants
and the embedding equality in Theorem \ref{discrete-embedding-rt}.

Finally, substituting estimates \refe{term-trunc-1},
\refe{term-nonlinear-1-final} and \refe{term-nonlinear-2-final} into
\refe{error-equation}, we obtain
\begin{align}
  & \left(D_\tau e_{u}^n , e_{u}^n \right) + \|e_{\bm{\sigma}}^n\|_{L^2}^2
  \nn \\
  & \le \epsilon \|e_{\bm{\sigma}}^{n-1}\|_{L^2}^2 + \epsilon^{-1} C \|e_{u}^{n-1}\|_{L^2}^2
  + \epsilon^{-1}C \|e_{u}^n\|_{L^2}^2 + \epsilon^{-1}C (\tau^2 + h^{2r+2}), 
\end{align}
Then, we chose a small $\epsilon$ and summing up the last inequality for
the index $n=1$, $2$, $\ldots$, $k$ to deduce that
\begin{align}
  & \| e_{u}^n \|_{L^2}^2 + \tau \sum_{m=1}^{n} \|e_{\bm{\sigma}}^m\|_{L^2}^2
  \le \tau C \sum_{m=1}^{n} \|e_{u}^m\|_{L^2}^2 
  + \tau C \sum_{m=1}^{n}(\tau^2 + h^{2r+2}), 
\end{align}
Thanks to the discrete Gronwall's inequality in Lemma \ref{gronwall},
when $C \tau \le \frac{1}{2}$, we have
\begin{align}
  & \| e_{u}^n \|_{L^2}^2 + \| e_{\bm{\sigma}}^n \|_{L^2}^2 
  + \tau \sum_{m=1}^{n} \left(\|e_{\bm{\sigma}}^m\|_{L^2}^2
  + \|\D e_{\bm{\sigma}}^m\|_{L^2}^2 \right)
  \nn \\
  & \le C\exp \left( \frac{TC}{1-C\tau} \right)
  (\tau^2+h^{2r+2})
  \nn \\
  & \le C \exp(2TC) (\tau^2+h^{2r+2})
\end{align}
Thus, \refe{mainresults} holds for $n=k$ if we take $\frac{C_*}{2} \ge C\exp(2TC)$.
We complete the induction. 

Theorem {\ref{maintheorem-artificial}} follows immediately from the the projection 
error estimates and the above inequality. \qed

\section{Numerical examples}
\label{sec-numer}
\setcounter{equation}{0}
In this section, we provide numerical 
experiments in both two and three dimensional spaces
to confirm our theoretical results in Theorem \ref{maintheorem-artificial}
and show the efficiency of the 
linearized mixed FEMs. 
The computations are carried out with
the free software FEniCS {\cite{fenics}}.

%================================================
%--- example  artificial, order, stability
%================================================
\begin{example}
\label{example-2d-square}
\rm
First we consider an artificial problem in two dimensional space 

\begin{empheq}[left=\empheqlbrace]{align} 
& \frac{\partial u}{\partial t} - \Delta u + u^3 = g\, ,
&& \mathrm{in}\, \Omega 
  \label{pdetest_2d} \\
& u = 0, && \mathrm{on}\, \partial \Omega 
  \label{pdetest_2d_bc} \\ 
& u = u_0(x), && \mathrm{in}\ \Omega ,
  \label{pdetest_2d_init}
\end{empheq}
where we take $\Omega = (0,1) \times (0,1)$.
The function $g$ is chosen 
correspondingly to the exact solution 
\begin{align}
  u(t,x,y) = \exp(t) x y (1-x) (1-y).
\nn
\end{align}
We set the terminal time $T = 1.0$ in this example.
This example had been tested in \cite{KimPB},
where a nonlinear backward Euler mixed FEM
with a two-grid algorithm was used.

\begin{figure}[htp]
\centering
\begin{tabular}{c}
\epsfig{file=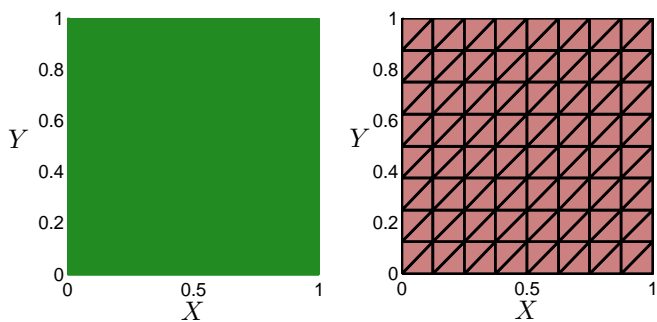,width=3.5in,height=1.7in}
\end{tabular}
\caption{A uniform triangular mesh on the unit square with $M=8$.}
\label{mesh2d}
\end{figure}
We use a uniform triangular mesh with $M+1$ vertices
in each direction, where $h = \frac{\sqrt{2}}{M}$
(see Figure \ref{mesh2d} for illustration with $M=8$).
We solve \refe{pdetest_2d}-\refe{pdetest_2d_init}
by the proposed linearized mixed FEMs
\refe{artificial-mixed-fem1}-\refe{artificial-mixed-fem2}
with $r=0$, $1$, $2$, respectively. 
To demonstrate the $O(\tau + h^{r})$ convergence of $L^{2}$-norm errors
of $u_h$ and $\bm{\sigma}_h$,
we set $\tau = \left( \frac{1}{M} \right)^{r+1}$ in our computation.
The $L^2$-norm errors of the scheme are shown in Table \ref{table_2d}.
From Table \ref{table_2d}, we can see that the $L^{2}$-norm errors  
for $u_h$ and $\bm{\sigma}_h$
are proportional to $h^{r+1}$,
which confirms the optimal convergence rates clearly.
{\setlength{\extrarowheight}{3pt}
\begin{table}[htp]
\begin{center}
\caption{$L^{2}$-norm errors of $u_h$ and $\bm{\sigma}_h$ on the unit square 
(Example \ref{example-2d-square}).}
\label{table_2d}
\begin{tabular}{c|c|c}
\hline
\hline
%===================================
%-  RT0 element 
%===================================
$\bH_h^{0} \times {V}_h^{0} 
\quad \tau = \frac{1}{M}$&
$ \|u^N-u_h^N\|_{L^2} $       &
$ \|\bm{\sigma}^N-\bm{\sigma}_h^N\|_{L^2} $  \\
\hline
$M=32 $ & 2.9850e-03 & 1.2659e-02   \\
$M=64 $ & 1.4928e-03 & 6.3329e-03   \\
$M=128$ & 7.4643e-04 & 3.1668e-03   \\
\hline
Order   & 9.9983e-01 & 9.9951e-01   \\
\hline
\hline
%===================================
%-  RT1 element 
%===================================
$\bH_h^{1} \times {V}_h^{1} 
\quad \tau = \frac{1}{M^2}$&
$ \|u^N-u_h^N\|_{L^2} $       &
$ \|\bm{\sigma}^N-\bm{\sigma}_h^N\|_{L^2} $       \\
\hline
$M=16$ & 2.3732e-04 & 1.0243e-03 \\
$M=32$ & 5.9385e-05 & 2.5731e-04 \\
$M=64$ & 1.4850e-05 & 6.4475e-05 \\
\hline                              
Order  & 1.9992e+00 & 1.9949e+00 \\
\hline
\hline
%===================================
%-  RT2 element 
%===================================
$\bH_h^{2} \times {V}_h^{2} 
\quad \tau = \frac{1}{M^3}$&
$ \|u^N-u_h^N\|_{L^2} $       &
$ \|\bm{\sigma}^N-\bm{\sigma}_h^N\|_{L^2} $       \\
\hline
$M=  8$ & 4.2973e-05 & 1.4866e-04    \\
$M= 16$ & 5.3949e-06 & 1.8728e-05    \\
$M= 32$ & 6.7509e-07 & 2.3501e-06    \\

\hline          
Order   & 2.9961e+00 & 2.9916e+00    \\
\hline
\hline
\end{tabular}
\end{center}
\end{table}
}

To test the stability of the proposed method, 
we solve \refe{pdetest_2d}-\refe{pdetest_2d_init} by the 
linearized mixed FEMs \refe{artificial-mixed-fem1}-\refe{artificial-mixed-fem1} 
with three fixed time steps $\tau=0.1$, $0.05$, $0.01$
on gradually refined meshes with  $M=8$, $16$, $32$, $64$ and $128$,
where we take $r=1$, i.e., $\bH_h^{1}(\Omega) \times {V}_h^{1}(\Omega)$ is used. 
We plot in Figure \ref{stab_d2} the $L^2$ errors of $u_h$ and $\bm{\sigma}_h$.
From Figure \ref{stab_d2}, we can see that for each fixed $\tau$, when  
the mesh is refined gradually, each $L^2$ error converges 
to a small constant of $O(\tau)$. This shows that the proposed 
linearized mixed FEM is unconditionally stable,
i.e., the method does not require mesh ratio restriction
$\tau \le C h^{\alpha}$ for a certain $\alpha>0$. 
\begin{figure}[htp]
\centering
\epsfig{file= 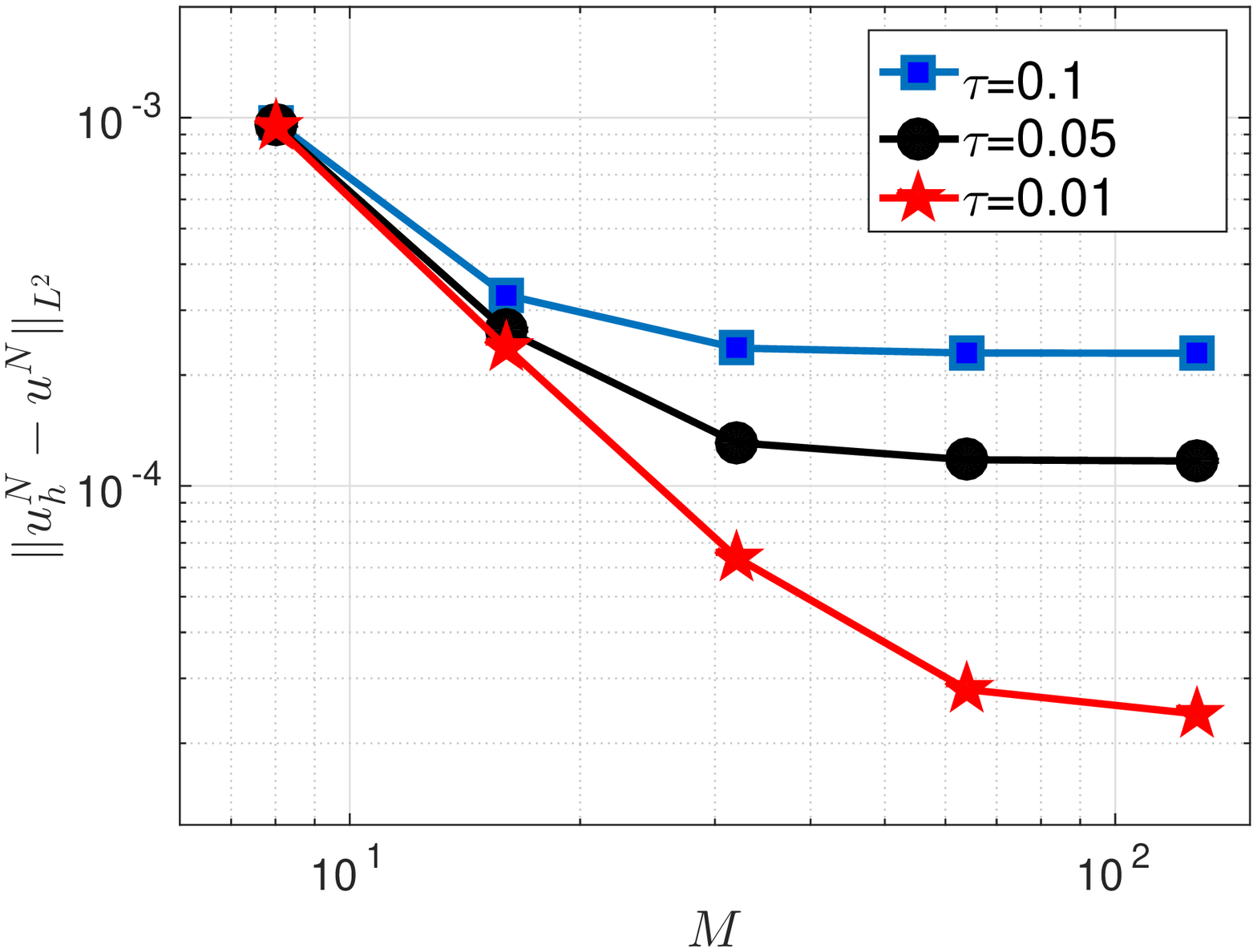,     width=1.8in,height=1.7in} 
\epsfig{file= 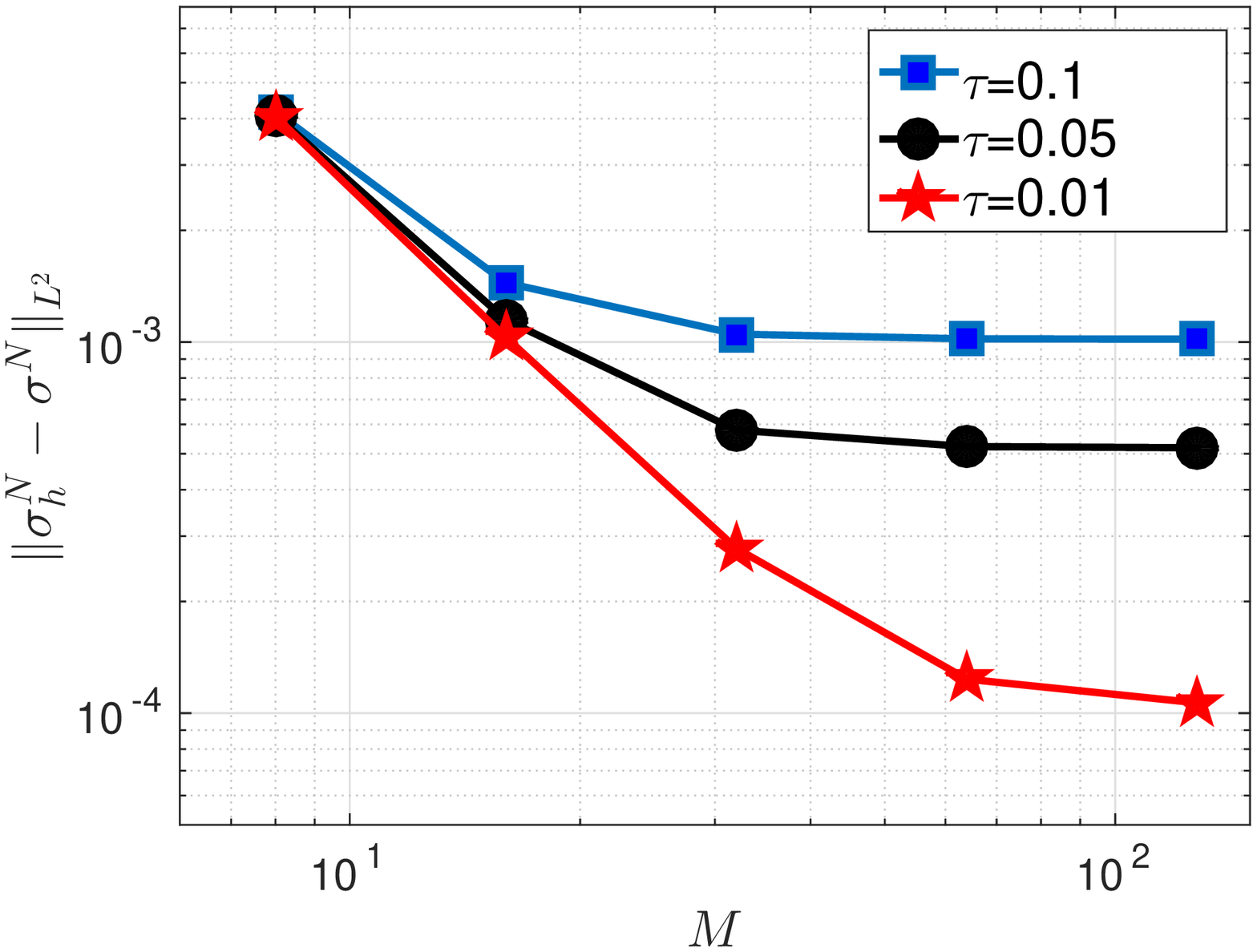,     width=1.8in,height=1.7in} 
\caption{$L^2$ errors of $u_h$ and $\bm{\sigma}_h$ with 
$\bH_h^{1} \times {V}_h^{1}$ on
gradually refined meshes with fixed $\tau$
(Example \ref{example-2d-square}).}
\label{stab_d2}
\end{figure}

\end{example}

%================================================
%--- example  artificial, order, stability 3d 3D
%================================================
\begin{example}
\label{example-3d-cube}
\rm
In this example we test the performance of the linearized mixed FEMs
for the following three-dimensional problem
\begin{empheq}[left=\empheqlbrace]{align} 
  & \frac{\partial u}{\partial t} - \Delta u + ({\bf{b}} \cdot \nabla u)u
  + u^3 - u = g\, ,
  && \mathrm{in}\, \Omega 
  \label{pdetest_3d} \\
  & u = 0, && \mathrm{on}\, \partial \Omega 
  \label{pdetest_3d_bc} \\ 
  & u = u_0(x), && \mathrm{in}\ \Omega ,
  \label{pdetest_3d_init}
\end{empheq}
where we take the unit cube $\Omega = (0,1) \times (0,1) \times (0,1)$.
The function $g$ is chosen 
correspondingly to the exact solution 
\begin{align}
  u(t,x,y) = \exp(-t) \sin(\pi x) \sin (2 \pi y) z(1-z).
\nn
\end{align}
A uniform tetrahedral mesh with $M+1$ vertices
in each direction are used,
where $h = \frac{\sqrt{3}}{M}$.
We solve the above equation \refe{pdetest_3d}-\refe{pdetest_3d_init}
by the proposed linearized mixed FEMs
\refe{artificial-mixed-fem1}-\refe{artificial-mixed-fem2}
with $(\bH_h^{r}(\Omega) \times {V}_h^{r}(\Omega) )$
for $r=0$, $1$, respectively. 
We also set $\tau = \left( \frac{1}{M} \right)^{r+1}$
and terminal time $T = 1.0$ in our computations.
We present the $L^2$-norm errors of the scheme in Table \ref{table_3d}.
Again, we can see clearly that the $L^{2}$-norm errors  
of $u_h$ and $\bm{\sigma}_h$
are proportional to $h^{r+1}$, for $r=0$, $1$, respectively.
This indicates that the convergence rate of the linearized mixed FEM
\refe{artificial-mixed-fem1}-\refe{artificial-mixed-fem2} is optimal
in three-dimensional space.
{\setlength{\extrarowheight}{3pt}
\begin{table}[htp]
\begin{center}
\caption{$L^{2}$-norm errors of $u_h$ and $\bm{\sigma}_h$ on the unit cube
(Example \ref{example-3d-cube}).}
\label{table_3d}
\begin{tabular}{c|c|c}
\hline
\hline
%===================================
%-  RT0 element 
%===================================
$\bH_h^{0} \times {V}_h^{0} 
\quad \tau = \frac{1}{M}$&
$ \|u^N-u_h^N\|_{L^2} $       &
$ \|\bm{\sigma}^N-\bm{\sigma}_h^N\|_{L^2} $  \\
\hline
$M=10 $ & 5.1823e-03  & 4.2003e-02     \\
$M=20 $ & 2.6285e-03  & 2.1121e-02     \\
$M=40 $ & 1.3189e-03  & 1.0575e-02     \\
\hline      
Order   & 9.8711e-01  & 9.9490e-01     \\
\hline
\hline
%===================================
%-  RT1 element 
%===================================
$\bH_h^{1} \times {V}_h^{1} 
\quad \tau = \frac{1}{M^2}$&
$ \|u^N-u_h^N\|_{L^2} $       &
$ \|\bm{\sigma}^N-\bm{\sigma}_h^N\|_{L^2} $       \\
\hline
$M= 8$ & 8.0631e-04  & 5.4993e-03   \\
$M=16$ & 2.0467e-04  & 1.3935e-03   \\
$M=32$ & 5.1364e-05  & 3.4997e-04   \\
\hline    
Order  & 1.9862e+00  & 1.9870e+00   \\
\hline
\hline
\end{tabular}
\end{center}
\end{table}
}

To test the stability of the proposed method, 
we solve \refe{pdetest_3d}-\refe{pdetest_3d_init} by the 
linearized mixed FEMs \refe{artificial-mixed-fem1}-\refe{artificial-mixed-fem1} 
with three fixed time steps $\tau=0.1$, $0.05$, $0.01$
on gradually refined meshes with  $M=10$, $20$, $30$, $40$ and $50$,
where $\bH_h^{1}(\Omega) \times {V}_h^{1}(\Omega)$ is used for spatial discretization. 
We plot in Figure \ref{stab_d3} the $L^2$ errors of $u_h$ and $\bm{\sigma}_h$.
From Figure \ref{stab_d3}, we can see that
the size of the time step $\tau$ affects the accuracy but not the stability
of the scheme.
This shows that the proposed 
linearized mixed FEMs are unconditionally stable in three-dimensional space.
\begin{figure}[htp]
\centering
\epsfig{file= 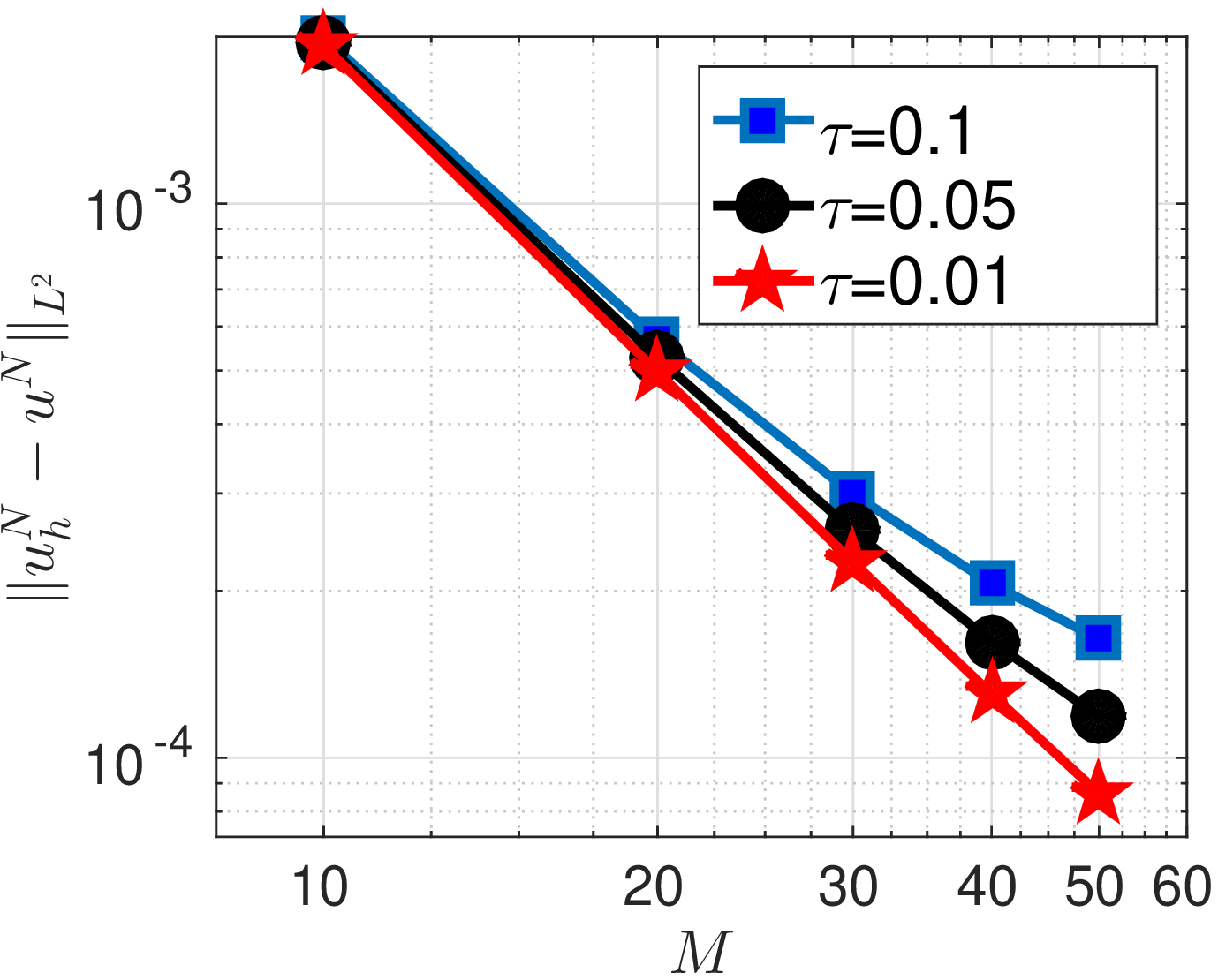,     width=1.8in,height=1.7in} 
\epsfig{file= 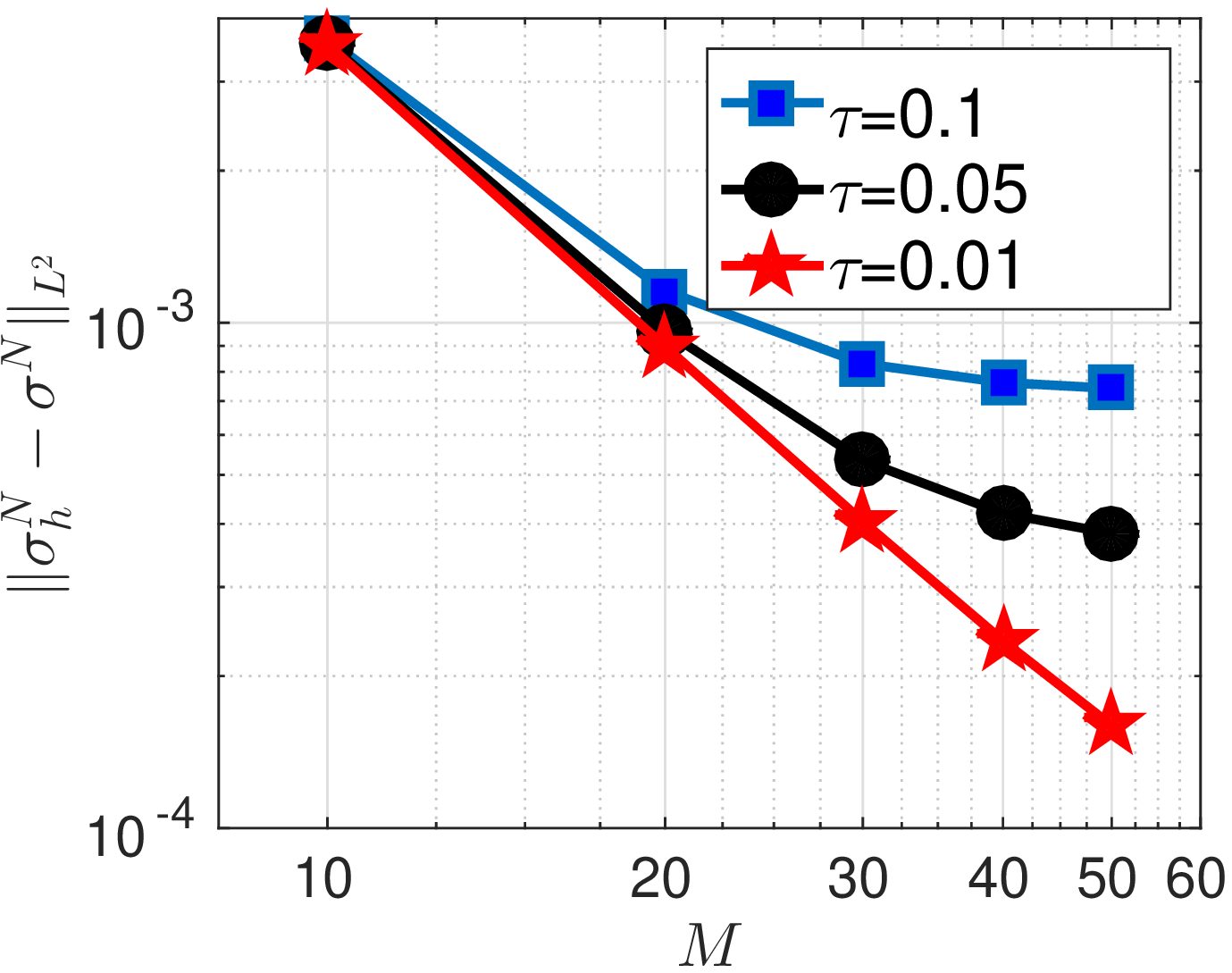,     width=1.8in,height=1.7in} 
\caption{$L^2$ errors of $u_h$ and $\bm{\sigma}_h$ with 
$\bH_h^{1} \times {V}_h^{1}$ on
gradually refined meshes with fixed $\tau$
(Example \ref{example-3d-cube}).}
\label{stab_d3}
\end{figure}
\end{example}

\section{Conclusion}

We have proved a discrete Sobolev embedding inequality
for the Raviart--Thomas mixed FEMs for
second order elliptic equations.
The essential idea is to control the $L^p$ norm
of $u_h$ by the discrete Sobolev norm $\|u_h\|_{DG}$
and then prove that $\|u_h\|_{DG}$ is bounded by $\|\bm{\sigma}_h\|_{L^2}$.
In this paper we focus on the Raviart--Thomas mixed FEMs.
However, it is easy to see that the results can be extended to
other stable elements,
such as Brezzi--Douglas--Marini(BDM) mixed FEMs.
We shall remark that in our proof there is no requirement on the domain $\Omega$.
In this paper, we only consider homogeneous Dirichlet boundary conditions.
It should be noted that extension to other boundary conditions can also be obtained
with slightly change of notations.
By using the proved discrete Sobolev inequality,
we have established an unconditionally optimal error estimates for mixed FEMs of
nonlinear parabolic equations.
We point out that the discrete Sobolev embedding inequalities proved in this work
can be used to analyze mixed FEMs of more general nonlinear parabolic systems.

\section*{Acknowledgments}
The authors would like to thank Prof.~Weiwei Sun
for useful discussions.

%==============================
%   references
%==============================

\end{document}